\documentstyle{article}
\title{A Note on the Cohomology of Moduli of Rank Two Stable Bundles}
\author{Richard Earl}
\date{March 12, 1996}
\begin{document}
\setcounter{section}{-1}
\newtheorem{lemma}{LEMMA}
\newtheorem{REM}[lemma]{Remark}
\newenvironment{rem}{\begin{REM} \rm}{\end{REM}}
\newtheorem{prop}[lemma]{PROPOSITION}
\newtheorem{cor}[lemma]{COROLLARY}
\newtheorem{thm}[lemma]{THEOREM}
\maketitle
\newcommand{\al}{\bar{\alpha}}
\newcommand{\ps}{\bar{\psi}}
\newcommand{\G}{{\cal G}}
\newcommand{\N}{{\cal N}}
\newcommand{\No}{{\cal N}_{0}}
\newcommand{\pa}{\overline{\partial}}
\newcommand{\g}{\bar{g}}
\newcommand{\Cmu}{{\cal C}_{\mu}}
\newcommand{\HG}{H^{*}_{\cal G}}
\newcommand{\Css}{{\cal C}^{s}}
\newcommand{\Q}{{\bf Q}}
\newcommand{\HS}{H^{*}}
\newcommand{\C}{{\cal C}}
\newcommand{\Cnu}{{\cal C}_{\nu}}
\newcommand{\GG}{\overline{\cal G}}
\newcommand{\HGG}{H^{*}_{\GG}}
\newcommand{\ax}{a_{1}^{1}}
\newcommand{\ay}{a_{1}^{2}}
\newcommand{\xii}{\tilde{\xi}}
\section{Introduction}
In recent years the cohomology ring of the moduli space ${\cal
N}_{g}(2,d)$ of rank two and odd degree $d$ stable bundles over a Riemann surface $M$ of genus $g \geq 2$ has been extensively studied \cite{B}, \cite{KN}, \cite{ST}, \cite{Z}. The subring $\HS({\cal N}_{g}(2,d);\Q)^{\Gamma}$ of $\HS({\cal N}_{g}(2,d);\Q)
$ which is invariant under the induced action of the mapping class group $\Gamma$ of $M$ has also been much studied and has been shown to play a central role in the ring structure of $\HS({\cal N}_{g}(2,d);\Q)$.\\
\indent In 1991 Zagier \cite{Z} began a study of certain relations in
the invariant cohomology ring. These are defined recursively in terms
of Newstead's generators $\alpha,\beta,\gamma$ \cite{N} by
\begin{equation}
(r+1)\zeta_{r+1} = \alpha \zeta_{r} + r \beta \zeta_{r-1} + 2 \gamma
\zeta_{r-2} \label{A}
\end{equation}
with $\zeta_{0}=1$ and $\zeta_{r}=0$ for $r<0$. Each of the authors
\cite{B}, \cite{KN}, \cite{ST}, \cite{Z} showed that $\zeta_{r}$ is a
relation in $\HS({\cal N}_{g}(2,d);\Q)^{\Gamma}$ for $r \geq g$ and
that $\zeta_{g}, \zeta_{g+1}, \zeta_{g+2}$  generate 
the relation ideal of the invariant cohomology ring. King and Newstead
further proved a decomposition theorem \cite[Prop.2.5]{KN}, originally
conjectured by Mumford, describing $\HS({\cal N}_{g}(2,d);\Q)$ in
terms of $\HS({\cal N}_{k}(2,d);\Q)^{\Gamma}$ $(k \leq g)$ and
exterior powers of $H^{3}({\cal N}_{g}(2,d);\Q).$\\   
\indent The methods employed by the authors \cite{B}, \cite{KN},
\cite{ST}, \cite{Z} differ greatly from Kirwan's original proof
\cite[$\S$2]{K} of Mumford's conjecture. Mumford introduced relations
in $\HS({\cal N}_{g}(2,d);\Q)$ which are constructed from the
vanishing Chern classes of a rank $2g-1$ bundle $\pi_{!}V$ over ${\cal
N}_{g}(2,d)$ and conjectured that these relations are complete. Zagier
\cite[$\S$6]{Z} showed that the relations $\zeta_{r}, r \geq g$ form a subset of the Mumford relations and for this reason we will refer to the vanishing $\zeta_{r}$ as the Zagier-Mumford relations.\\
\indent The purpose of this note is two-fold. We will firstly
rederive the result that the first three Zagier-Mumford relations form
a minimal complete set for the invariant cohomology. The second result is to
prove a subsequent and stronger version of Mumford's conjecture;
namely we will show that the relations constructed solely from the
first vanishing Chern class $c_{2g}(\pi_{!}V)$ freely generate the
relation ideal of $\HS({\cal N}_{g}(2,d);\Q)$ as a $\Q[\alpha,\beta]$-module. Both results follow easily
from Kirwan's calculations in \cite[$\S$2]{K}. Partly the aim of this
note is to demonstrate the power of the methods of \cite{K} which
currently is the only approach to have generalised to the rank three
case \cite{E}.\\
\indent For ease of notation we will from now on write ${\cal N}$ for
${\cal N}_{g}(2,d)$ and write ${\cal N}_{0}$ for the moduli space of rank
two odd degree stable bundles of fixed determinant. Also we write $\g$
for $g-1$ and $[2g]$ for the set $\{1,...,2g\}$.
\section{Kirwan's Approach}
\indent Let $\C$ denote the space of all holomorphic structures on a
fixed $C^{\infty}$ complex vector bundle ${\cal E}$ over $M$ of rank two and
odd degree $d$ and let $\G_{c}$ denote the group of all $C^{\infty}$
complex automorphisms of ${\cal E}$. We may then identify $\N$ with
the quotient $\Css/\G_{c}$ where $\Css \subset \C$ is the open subset
consisting of stable holomorphic structures. Let $\G$ denote the gauge group of
all $C^{\infty}$ automorphisms of ${\cal E}$ which are unitary with
respect to a fixed Hermitian structure and let $\overline{\G}$ denote
the quotient of $\G$ by its $U(1)$ centre.\\
\indent Then $\HS(\N;\Q)$ is naturally isomorphic to
$H^{*}_{\overline{\G}}(\Css;\Q)$ \cite[9.1]{AB}, and Atiyah and Bott
show further \cite[thm.7.14]{AB} that the restriction map
\begin{equation}
\HS(B\GG;\Q) \cong H^{*}_{\overline{\G}}(\C;\Q) \rightarrow
H^{*}_{\overline{\G}}(\Css;\Q) \cong \HS(\N;\Q) \label{B}
\end{equation}
is surjective. They construct a rank two $\G$-equivariant
holomorphic bundle ${\cal V}$ over $\C \times M$ and define generators
\begin{equation}
a_{1}, a_{2}, f_{2} \mbox{ and } b_{1}^{s}, b_{2}^{s} \quad (s \in
[2g]) \label{gen}
\end{equation}
for $\HG(\C;\Q) = \HS(B\G;\Q)$ by
taking the slant products
\[
c_{r}({\cal V}) = a_{r} \otimes 1 + \sum_{s=1}^{2g} b_{r}^{s} \otimes
e_{s} + f_{r} \otimes \omega
\]
where $e_{1},...,e_{2g}$ is a fixed basis for $H^{1}(M;\Q)$ and
$\omega$ is the standard generator for $H^{2}(M;\Q)$. The only
relations amongst the generators (\ref{gen}) are that $a_{1}, a_{2},
f_{2}$ commute with everything and that the $b_{i}^{s}$ anticommute
amongst themselves.\\
\indent Rather than (\ref{B}) we shall consider the surjection
\begin{equation}
H^{*}_{\overline{\G}}(\C;\Q) \otimes \Q[a_{1}] \cong \HG(\C;\Q) \rightarrow
\HG(\Css;\Q) \cong \HS(\N;\Q) \otimes \Q[a_{1}]. \label{C}
\end{equation}
The images of (\ref{gen}) under this map form generators for
$\HG(\Css;\Q)$ which we will also refer to as 
$a_{1}, a_{2}, f_{2}, b_{1}^{s},$ $b_{2}^{s}$ and the relations among
these restrictions form the kernel of (\ref{C}). In order to study
this kernel we introduce a $\G$-perfect stratification of $\C$ due to
Shatz \cite{Sh}.\\
\indent 
\indent Any unstable holomorphic bundle $E$ over $M$ of rank $n$ and degree $d$
has a canonical filtration (or flag) \cite[p.221]{HN} \label{filt}
which in the rank $n=2$ case is a line subbundle $L$ of $E$ of
degree $d_{1}$ such that $d_{1} > d/2$. We define the type of $E$ to be $(d_{1},d-d_{1})$ and define the
type of a stable bundle to be $\mu_{0} = (d/2,d/2)$. The stratum $\Cmu
\subseteq {\cal C}$ is the set of all 
holomorphic vector bundles of type $\mu$ and we construct a total
order $\preceq$ on the set of types by writing $(\mu_{1},\mu_{2}) \preceq
(\nu_{1},\nu_{2})$ if $\mu_{1} \leq \nu_{1}$.\\ 
\indent Kirwan's proof of Mumford's conjecture is based upon a set of
completeness criteria for a set $\cal R$ of relations in $\HG(\C;\Q)$
\cite[Prop.1]{K}. These criteria involve finding for each $\mu \neq
\mu_{0}$ relations ${\cal R}_{\mu} \subseteq {\cal R}$ which in a
technical sense correspond to the stratum $\Cmu$. We introduce here
similar completeness criteria for the invariant cohomology:
\begin{prop}
\label{ICC}
(Invariant Completeness Criteria)  Let ${\cal R}$
be a subset of the kernel of the restriction map
\begin{equation}
\HG(\C;\Q)^{\Gamma} \rightarrow \HG(\Css;\Q). \label{Q23}
\end{equation}
Suppose that for each unstable type $\mu$ there is a subset ${\cal
R}_{\mu}$ of the ideal generated by ${\cal R}$ in $\HG(\C;\Q)^{\Gamma}$ such
that the restriction of ${\cal R}_{\mu}$ to $\HG(\C_{\nu};\Q)$ 
is zero when $\nu \prec \mu$ and when $\nu = \mu$ equals the ideal
generated by $e_{\mu}$ in $\HG(\C_{\mu};\Q)^{\Gamma}$, where $e_{\mu}$
denotes the equivariant Euler class of the normal bundle to $\Cmu$ in
$\C$. Then ${\cal R}$ generates the kernel of the restriction map
(\ref{Q23}) as an ideal of $\HG(\C;\Q)^{\Gamma}.$ 
\end{prop}
{\bf PROOF:} We include now the main points in the proof of the above
proposition. However the only difference between this proof and the
argument of \cite[prop.4]{E} is to observe that $\Cmu$ is
$\Gamma$-invariant and hence $e_{\mu} \in \HG(\Cmu;\Q)^{\Gamma}$.\\
\indent For $\mu$ an unstable type let $\mu-1$ denote
the type previous to $\mu$ with respect to $\preceq$ and define
$V_{\mu} = \bigcup_{\nu \preceq \mu} \Cnu.$ Then $V_{\mu}$ is an open
subset of $\C$ which contains $\Cmu$ as a closed submanifold.\\
\indent Let $d_{\nu}$ denote the complex codimension of $\Cnu$ in
$\C$. For any given $i \geq 0$ there are only finitely many $\nu \in
{\cal M}$ such that $2d_{\nu} \leq i$ \cite[7.16]{AB} and so for each $i \geq 0$ there
exists some $\mu$ such that 
\[
H_{\G}^{i}(\C;\Q) = H_{\G}^{i}(V_{\mu};\Q).
\]
Hence it is enough to show that for each $\mu$ the image in
$\HG(V_{\mu};\Q)^{\Gamma}$ of the ideal generated by ${\cal R}$ contains the
image in $\HG(V_{\mu};\Q)^{\Gamma}$ of the kernel of (\ref{Q23}). Note that the
above is clearly true for $\mu=\mu_{0}$ as $V_{\mu_{0}}=\Css.$ We will
proceed by induction with respect to $\preceq$.\\
\indent Assume now that $\mu \neq \mu_{0}$ and that $\zeta \in \HG(\C;\Q)^{\Gamma}$ lies in the kernel of (\ref{Q23}). Suppose that the image of $\zeta$ in $\HG(V_{\mu-1};\Q)$ is in the
image of the ideal generated by ${\cal R}.$ We may, without any loss of
generality, assume that the image of $\zeta$ in $\HG(V_{\mu-1};\Q)$ is
zero. Then by the exactness of the Thom-Gysin sequence
\[
\cdot \cdot \cdot \rightarrow H_{\G}^{*-2d_{\mu}}(\C_{\mu};\Q)
\rightarrow \HG(V_{\mu};\Q)
\rightarrow \HG(V_{\mu-1};\Q) \rightarrow \cdot \cdot \cdot
\]
there exists an
element $\eta \in H_{\G}^{*-2d_{\mu}}(\C_{\mu};\Q)$ which is
mapped to the image of $\zeta$ in $\HG(V_{\mu};\Q)$ by the Thom-Gysin
map. The composition 
\[
H_{\G}^{*-2d_{\mu}}(\C_{\mu};\Q) \rightarrow \HG(V_{\mu};\Q)
\rightarrow \HG(\C_{\mu};\Q) 
\]
is given by multiplication by $e_{\mu}$ which is
not a zero-divisor in $\HG(\C_{\mu};\Q)$ \cite[p.569]{AB}. Hence the
restriction of $\zeta$ in $\HG(\C_{\mu};\Q)$ is $\eta e_{\mu}$ and by
our initial observation $\eta 
\in H^{*-2d_{\mu}}_{\G}(\Cmu;\Q)^{\Gamma}.$ By hypothesis there exists
$\theta$ in ${\cal R}_{\mu}$ whose image in $\HG(\Cmu;\Q)$
equals $\eta e_{\mu}.$\\
\indent Kirwan shows \cite[p.867]{K} that the direct sum of
restriction maps 
\begin{equation}
\HG(V_{\mu};\Q) \rightarrow \bigoplus_{\nu \preceq \mu}
\HG(\C_{\nu};\Q) \label{sumres}
\end{equation}
is injective. The images of $\theta$ and $\zeta$ under (\ref{sumres})
are equal and hence the restrictions of $\theta$ and $\zeta$ to 
$\HG(V_{\mu};\Q)$ are the same, completing the proof. $\indent \Box$
\section{The Mumford and Zagier-Mumford Relations}
\indent The group $\GG$ acts freely on $\Css$ and the $U(1)$-centre of
$\G$ acts as scalar multiplication on the fibres of ${\cal V}$. The projective
bundle of ${\cal V}$ descends to a holomorphic projective bundle over
${\cal N} \times M$ which is the projective bundle of a universal
holomorphic bundle $V$ of rank two and odd degree \cite[p.580]{AB}.\\
\indent The bundle $V$ is universal in the sense that the restriction
of $V$ to $[E] \times M$ for any class $[E] \in {\cal N}$ is
isomorphic to $E$. Note V is not unique; we may tensor $V$ by the pullback
of a line bundle over ${\cal N}$ to produce a second bundle with the same
universal property. However we may normalise $V$ \cite[p.857, p.877]{K} by
requiring the relation
\begin{equation}
f_{2} = (d - 2 \g)a_{1} + \sum_{s=1}^{g}
b_{1}^{s} b_{1}^{s+g}. \label{norm}
\end{equation}
\indent Let $\pi: \N \times M \to \N$ be the first projection. When $d
= 4g-3$ then any $E \in \Css$ has slope $\mu(E) = d/n > 2\g$ and thus
\cite[lemma 5.2]{N2} $H^{1}(M,E) = 0$. Hence $\pi_{!}V$ is a genuine vector
bundle over $\N$ of rank $2g-1$ with fibre $H^{0}(M,E)$ over $[E] \in
\N$. We know from \cite[prop.9.7]{AB} that
\[
\HG(\Css;\Q) \cong \HS({\cal N}_{0};\Q) \otimes \Q [a_{1}] \otimes
\Lambda^{*}\{b_{1}^{1},...,b_{1}^{2g}\}.
\]
The Mumford relations $c_{r,S}$ $(r \geq 2g, S \subseteq [2g])$ are
then defined by writing
\begin{equation}
c_{r}(\pi_{!}V) = \sum_{S \subseteq [2g]} c_{r,S} \prod_{s \in S}
b_{1}^{s}, \label{first}
\end{equation}
where each $c_{r,S}$ is written in terms of generators for $\HS({\cal
N}_{0};\Q) \otimes \Q [a_{1}]$, namely $a_{1}$ and Newstead's
generators $\alpha, \beta, \psi_{s}$. In terms of
the generators (\ref{gen}) these are given by
\[
\alpha = 2f_{2} - da_{1},\indent \beta = (a_{1})^{2} - 4a_{2}, \indent
\psi_{s} = 2b_{2}^{s}.
\]
Kirwan's proof of Mumford's conjecture \cite[$\S$2]{K} shows that the
Mumford relations together with the normalising relation (\ref{norm})
form a complete set of relations for $\HS({\cal N}_{0};\Q)$. Following
Kirwan \cite[p.871]{K} we reformulate the definition (\ref{first}) and
write
\[
\Psi(t) = \sum_{r=0}^{\infty} c_{r}(\pi_{V}) t^{2g-1-r} = \sum_{r=-\infty}^{\g}
(\sigma_{r}^{0} + \sigma_{r}^{1} t) (t^{2} + a_{1} t + a_{2})^{r},
\indent \sigma_{r}^{k} = \sum_{S \subseteq [2g]} \sigma_{r,S}^{k}
\prod_{s \in S} b_{1}^{s}.
\]
We will also refer to $\sigma_{r,S}^{k} (k=0,1, r<0, S \subseteq
[2g])$ as the Mumford relations. (Note $\sigma_{r}^{0}$ and
$\sigma_{r}^{1}$ differ slightly from Kirwan's terms $\sigma_{r}$ and
$\tau_{r}$.) This new formulation will prove more
convenient when we need to determine the restrictions of the Mumford
relations to various strata. Atiyah and Bott define generators 
\[
a_{1}^{1}, a_{1}^{2}, \mbox{ and } b_{1}^{1,s}, b_{1}^{2,s}
\]
for $\HG(\Cmu;\Q)$ via the isomorphism \cite[prop.7.12]{AB}
\[
\HG(\Cmu;\Q) \cong H^{*}_{\G(1,d_{1})}(\C(1,d_{1})^{ss};\Q) \otimes
H^{*}_{\G(1,d_{2})}(\C(1,d_{2})^{ss};\Q).
\]
In terms of these generators the crucial calculation of Kirwan in her
proof of Mumford's conjecture is:
\begin{lemma}[Kirwan]
\rm{\cite[pp.871-873]{K}} {\em Let $\mu = (d_{1},d_{2})$ and write $D = d_{2} - 2g +1$. Then the
restrictions $\sigma_{D,S}^{k,\mu}$ of $\sigma_{D,S}^{k}$ $(k=0,1)$ in
$\HG(\Cmu;\Q)$ are given by
\begin{equation}
\sigma_{D,S}^{0,\mu}  = \frac{(-1)^{g\g/2}}{2^{2g} g!} \left( \prod_{s
\not \in S} (b_{1}^{2,s} - b_{1}^{1,s}) \right) a_{1}^{1} e_{\mu}, \indent
\sigma_{D,S}^{1,\mu} = \frac{(-1)^{g\g/2}}{2^{2g} g!} \left( \prod_{s
\not \in S} (b_{1}^{2,s} - b_{1}^{1,s}) \right) e_{\mu}. \label{Q61}
\end{equation}}
\end{lemma}
This calculation plays a major role in the following two theorems.
\begin{thm}
(\cite{B}, \cite{KN}, \cite{ST}, \cite{Z}.)
Each of the sets
\begin{eqnarray}
\{ \zeta_{g}, \zeta_{g+1}, \zeta_{g+2} \}, \label{Q43}\\
\{ c_{2g,[2g]}, c_{2g+1,[2g]}, c_{2g+2,[2g]} \}, \label{Q44} \\
\{ \sigma_{-1,[2g]}^{1}, \sigma_{-1,[2g]}^{0}, \sigma_{-2,[2g]}^{1}
\}, \label{Q45}
\end{eqnarray} 
forms a minimal complete set of relations for the invariant cohomology
ring $H^{*}({\cal N}_{0};\Q)^{\Gamma}$.
\end{thm}
{\bf PROOF:} From \cite[$\S$6]{Z} we know that the relations (\ref{Q43}) and
(\ref{Q44}) generate the same ideal of $\HG(\C;\Q)^{\Gamma}$. Since
\begin{eqnarray*}
c_{2g,[2g]} & = & \sigma_{-1,[2g]}^{1},\\
c_{2g+1,[2g]} & = & \sigma_{-1,[2g]}^{0} -a_{1} \sigma_{-1,[2g]}^{1},\\
c_{2g+2,[2g]} & = & \sigma_{-2,[2g]}^{1} -a_{1} \sigma_{-1,[2g]}^{0}
+((a_{1})^{2} -a_{2}) \sigma_{-1,[2g]}^{1},
\end{eqnarray*}
we can see that the relations (\ref{Q45}) also generate the same
ideal.\\
\indent Let $\mu = (d_{1},d_{2})$. Now $\HG(\Cmu;\Q)^{\Gamma}$ is generated by
\[
\ax, \ay, \xi_{1,1}^{1,1}, \xi_{1,1}^{1,2}+ \xi_{1,1}^{2,1},
\xi_{1,1}^{2,2},
\]
where $\xi_{1,1}^{i,j} = \sum_{s=1}^{g} b_{1}^{i,s}
b_{1}^{j,s+g}$. On restriction to $\HG(\Cmu;\Q)$
\[
a_{1} \mapsto a_{1}^{1} + a_{1}^{2}, \qquad b_{1}^{s} \mapsto b_{1}^{1,s}
+ b_{1}^{2,s}, \qquad f_{2} \mapsto d_{1}a_{1}^{2} + d_{2}a_{1}^{1} +
\xi_{1,1}^{1,2} + \xi_{1,1}^{2,1}.
\]
Hence $e_{\mu}\HG(\Cmu;\Q)^{\Gamma}$ is generated by
\begin{equation}
\sum_{ \begin{array}{c} \mbox{\scriptsize $S \subseteq [g]$} \\[-2pt]
\mbox{ \scriptsize $|S| = k$} \end{array}} \hspace{-4mm}
 \left( \prod_{s \in S}
b_{1}^{1,s}b_{1}^{1,s+g}-b_{1}^{2,s}b_{1}^{2,s+g} \right) (a_{1}^{1})^{i}
e_{\mu} \label{mess}
\end{equation}
for $i = 0,1,$ $0 \leq k \leq g$ and the restrictions of $a_{1}$,
$f_{2}$ and $\xi_{1,1} = \sum_{s=1}^{g} b_{1}^{s} b_{1}^{s+g}$. Let $P(S)$
denote the set of partitions $S$ into two sets $S_{1},S_{2}$. We then
see from lemma 2 that (\ref{mess}) above is the restriction of
\begin{equation}
\frac{1}{2^{k}} \hspace{-4mm} \sum_{ \begin{array}{c}
\mbox{\scriptsize $S \subseteq  [g]$} \\[-2pt] 
\mbox{ \scriptsize $|S| = k$} \end{array}} \hspace{-3mm}
\sum_{P(S)} (\pm)
\left( \prod_{s \in S_{1}} b_{1}^{s} \right)
\left( \prod_{s \in S_{2}} b_{1}^{s+g} \right)
\sigma^{i}_{D,[2g] - (S_{2} \cup (S_{1}+g))} \label{mess2}
\end{equation}
where the sign $(\pm)$ depends on the particular partition of
$S$. Thus by proposition \ref{ICC} the relations (\ref{mess2}) above
for $D < 0,$ $i=0,1,$ $0 \leq k \leq g$ generate the invariant relation ideal of
\[
\HG(\Css;\Q) \cong \HS({\cal N}_{0};\Q) \otimes
\Lambda^{*}\{b_{1}^{1}, \ldots ,b_{1}^{2g}\} \otimes \Q[a_{1}].
\]
In particular the relations
\[
\left\{ \sigma_{D,[2g]}^{i} : D<0,i=0,1 \right\}
\]
generate $\HS({\cal N}_{0};\Q)^{\Gamma}$. It follows from (\ref{A})
that the sets (\ref{Q43}), 
(\ref{Q44}) and (\ref{Q45}) each form a complete set of relations for
the invariant cohomology ring of ${\cal N}_{0}$.\\
\indent Minimality then follows easily. Suppose that for some $\eta,
\theta \in \HS({\cal N}_{0};\Q)^{\Gamma}$ we have
\begin{equation}
\sigma_{-2,[2g]}^{1} + \eta \sigma_{-1,[2g]}^{0} + \theta
\sigma_{-1,[2g]}^{1} = 0. \label{Q58}
\end{equation}
So $\eta$ has degree 2 and $\theta$ has degree 4.\\
\indent Let $\mu = (2\g+1,2\g)$. Restricting equation (\ref{Q58}) to
$\HG(\Cmu;\Q)$ we find from (\ref{Q61}) that
\begin{equation}
\eta_{\mu} a_{1}^{1} + \theta_{\mu} = 0 \label{Q59}
\end{equation}
since $e_{\mu}$ is not a zero-divisor in $\HG(\Cmu;\Q)$. The
restriction map
\begin{equation}
\HG(\C;\Q) \rightarrow \HG(\Cmu;\Q) \label{Q60}
\end{equation}
is given by
\[
a_{1} \mapsto a_{1}^{1} + a_{1}^{2}, \indent a_{2} \mapsto a_{1}^{1}
a_{1}^{2}, \indent f_{2} \mapsto 2\g a_{1}^{1} + (2\g+1) a_{1}^{2}
+ \xi_{1,1}^{1,2} + \xi_{1,1}^{2,1},
\]
\[
b_{1}^{s} \mapsto b_{1}^{1,s} + b_{1}^{2,s}, \indent b_{2}^{s} \mapsto
a_{1}^{2} b_{1}^{1,s} + a_{1}^{1} b_{1}^{2,s}.
\]
From (\ref{Q59}) we see that $\eta_{\mu}$ and $\theta_{\mu}$ are both
zero. Since the restriction map (\ref{Q60}) is injective in degrees 4
and less, we have that $\eta$ and $\theta$ are both zero -- which
contradicts (\ref{Q58}).\\
\indent Similarly the equation
\[
\sigma_{-1,[2g]}^{0} + \eta \sigma_{-1,[2g]}^{1} = 0
\]
has no solutions for $\eta \in H^{2}_{\G}(\C;\Q)$. \indent $\Box$
\section{Mumford's Conjecture}
The final result of this note is a stronger version of Mumford's
conjecture as proven by Kirwan \cite[$\S$2]{K} and which is
confusingly also referred to as Mumford's conjecture.\\
\indent The Poincar\'{e} polynomial of the relation ideal of
$\HS({\cal N}_{0};\Q)$ equals \cite[p.593]{AB}
\[
\frac{t^{2g}(1+t)^{2g}}{(1-t^{2})(1-t^{4})}.
\]
Now there are ${2g \choose r}$ relations of the form $c_{2g,S}$ of degree
$2g+r$; $\alpha$ has degree two, $\beta$ has degree four and neither
are nilpotent in $\HG(\C;\Q)$. This strongly suggests:
\begin{thm}
The relation ideal of $\HS({\cal N}_{0};\Q)$ is freely generated as a
$\Q[\alpha,\beta]$-module by the Mumford relations $c_{2g,S}$ for $S
\subseteq [2g]$.
\end{thm}
{\bf PROOF:} Define
\[
\al = \alpha - \sum_{s=1}^{g} b_{1}^{s} b_{1}^{s+g} = 2f_{2} - da_{1}
- \sum_{s=1}^{g} b_{1}^{s} b_{1}^{s+g}.
\]
We will show that the relations $c_{2g,S}$ generate the relation ideal
of $\HG(\Css;\Q)$ as a $\Q[\al,\beta]$-module. As $\al$ restricts to
$\alpha$ in $H^{*}({\cal N}_{0};\Q)$ then this is equivalent to
the above result. It will be sufficient to prove that
\begin{equation}
\sum_{S \subseteq [2g]} \lambda_{S}(\al,\beta) c_{2g,S} = 0
\indent \lambda_{S}(\al,\beta) \in \Q[\al,\beta] \label{Q53}
\end{equation}
in $\HG(\C;\Q)$ if and only if $\lambda_{S}(\al,\beta) = 0$
for each $S \subseteq [2g].$\\
\indent Let $\mu = (2\g+1,2\g)$. Then from lemma 2 we know that the
restriction of $c_{2g,S} = \sigma_{-1,S}^{1}$ in $\HG(\Cmu;\Q)$ equals
\[
\frac{(-1)^{g\g/2}}{2^{2g} g!} \left( \prod_{s
\not \in S} (b_{1}^{2,s} - b_{1}^{1,s}) \right) e_{\mu}.
\]
If we restrict equation (\ref{Q53}) to $\HG(\Cmu;\Q)$ and recall that
$e_{\mu}$ is not a zero-divisor in $\HG(\Cmu;\Q)$ \cite[p.569]{AB} we obtain 
\[
\sum_{S \subseteq [2g]} \lambda_{S}(\al_{\mu},\beta_{\mu}) \left( \prod_{s \not
\in S} (b_{1}^{2,s}-b_{1}^{1,s}) \right) = 0.
\]
Now the restrictions of $\al$ and $\beta$ in $\HG(\Cmu;\Q)$ equal
\begin{equation}
\al_{\mu} = (a_{1}^{2}-a_{1}^{1}) - \sum_{s=1}^{g}
(b_{1}^{1,s} - b_{1}^{2,s})(b_{1}^{1,s+g}-b_{1}^{2,s+g}), \indent
\beta_{\mu} = (a_{1}^{2} - a_{1}^{1})^{2}. \label{Q54}
\end{equation}
By comparing the coefficients of $\prod_{s \in S} (b_{1}^{2,s} -
b_{1}^{1,s})$ for each $S \subseteq [2g]$ we see that
\[
\lambda_{S}(\al_{\mu},\beta_{\mu}) = 0 \indent S \subseteq [2g].
\]
\indent Consider the restriction map
\begin{equation}
\Q[\al,\beta] \rightarrow \Q[\al_{\mu},\beta_{\mu}]. \label{Q55}
\end{equation}
From the expressions (\ref{Q54}) of $\al_{\mu}$ and $\beta_{\mu}$
we can see that the kernel of the restriction map (\ref{Q55}) is
the ideal of $\Q[\al,\beta]$ generated by
\[
(\al^{2} - \beta)^{g+1}.
\]
However $\al^{2} - \beta$ is not a zero-divisor in
$\HG(\C;\Q)$. So we can assume without any loss of generality that for
some $S,$ $\lambda_{S}$ is either zero or not in the ideal generated by
$\al^{2}-\beta$. For this $S$ we have $\lambda_{S} = 0$ since
$\lambda_{S}(\al_{\mu},\beta_{\mu}) = 0$. Inductively we can see that
\[
\lambda_{S}(\al,\beta) = 0 \mbox{ for } S \subseteq [2g]. \indent \Box
\]   

Mathematical Institute,\\
24-29 St. Giles,\\
Oxford,\\
OX1 3LB,\\
England,\\[\baselineskip]
earl@maths.ox.ac.uk

\begin{thebibliography}{WW}
\bibitem{AB} M.F.Atiyah and R.Bott {\em The Yang-Mills equations over
Riemann surfaces} Philos. Trans. Roy. Soc. London Ser. A {\bf 308}
(1982) 523-615.
\bibitem{B} V.Baranovsky {\em Cohomology ring of the moduli space of
stable vector bundles with odd determinant} Izv. Russ. Acad. Nauk.
{\bf 58 n4} (1994) 204-210.
\bibitem{E} R.A.Earl {\em The Mumford relations and the moduli of rank three stable bundles} Compositio Math. (to appear).
\bibitem{HN} G.Harder and M.S.Narasimhan {\em On the cohomology groups
of moduli spaces of vector bundles over curves} Math. Ann. {\bf 212} (1975)
215-248.
\bibitem{KN} A.D.King and P.E.Newstead {\em On the cohomology of the
moduli space of rank 2 vector bundles on a curve} (1995 preprint).
\bibitem{K} F.C.Kirwan {\em Cohomology rings of moduli spaces of
bundles over Riemann surfaces} J. Amer. Math. Soc. {\bf 5} (1992) 853-906.
\bibitem{N} P.E.Newstead {\em Characteristic classes of stable
bundles of rank 2 over an algebraic curve} Trans. Amer. Math. Soc.
{\bf 169} (1972) 337-345.
\bibitem{N2} P.E.Newstead {\em Introduction to moduli problems and
orbit spaces} Tata Inst. Lect. 51 1978.
\bibitem{Sh} S.S.Shatz {\em The decomposition and specialization of
algebraic families of vector bundles} Compositio Math. {\bf 35}
(1977) 163-187.
\bibitem{ST} B.Siebert and G.Tian {\em Recursive relations for the
cohomology ring of moduli spaces of stable bundles} (1994 preprint).
\bibitem{Z} D.Zagier {\em On the cohomology of moduli spaces of rank 2
vector bundles over curves} (1995 preprint).
\end{thebibliography}
\end{document}